\documentclass[prl,twoside,onecolumn]{revtex4}
\usepackage{color}
\usepackage{graphicx} 

\begin{document}
\title{Ropelength and writhe quantization of 12-crossing knots}

\author{Alexander R. Klotz }
\author{Caleb J. Anderson}
\affiliation{Department of Physics and Astronomy, California State University, Long Beach}
\begin{abstract}
The ropelength of a knot is the minimum length required to tie it. Computational upper bounds have previously been computed for every prime knot with up to 11 crossings. Here, we present ropelength measurements for the 2176 knots with 12 crossings, of which 1288 are alternating and 888 are non-alternating. We report on the distribution of ropelengths within and between crossing numbers, as well as the space writhe of the tight knot configurations. It was previously established that tight alternating knots have a ``quantized'' space writhe close to a multiple of 4/7. Our data supports this for 12-crossing alternating knots and we find that non-alternating knots also show evidence of writhe quantization, falling near integer or half-integer multiples of 4/3, depending on the parity of the crossing number. Finally, we examine correlations between geometric properties and topological invariants of tight knots, finding that the ropelength is positively correlated with hyperbolic volume and its correlates, and that the space writhe is positively correlated with the Rasmussen s invariant.
\end{abstract}

\maketitle

\section{Introduction}
An ideal knot is the tightest possible configuration of a given knot. More formally, if a knot is tied in an incompressible, inextensible unit-radius tube with the knot along its axis, the ideal configuration of the tube has the minimum possible ratio of contour length to tube radius that respects a no-overlap constraint \cite{cantarella2002minimum, pieranski1998search}. The contour length of an ideal knot is known as the ropelength, which is a geometric invariant of the knot. There are existing upper and lower bounds on the scaling of ropelength with respect to the crossing number of the knot, \cite{diao2020braid, buck1998four, diao2019ropelengths}. Ropelength upper bounds of specific knots can be calculated numerically by initializing a discrete knot and perturbing its coordinates in a systematic way towards the minimum-length configuration. Several groups have exhaustively computed the ropelength of every prime knot with up to 10 crossings \cite{pieranski1998search, ashton2011knot}, and Brian Gilbert, a contributor to the Knot Atlas wiki, has computed the ropelength of all knots and links with up 11 crossings \cite{katlas}. In a previous work, we computed ropelengths for selected satellite and torus knots up to 1023 crossings, and we direct readers to the introduction of that work and references therein for a summary of ropelength bounds and measurements \cite{klotz2021ropelength}. Since its publication, it has been proven by Yuanan Diao that alternating knots must grow at least linearly with crossing number \cite{diao2022ropelength}, and a linear upper bound has been established for 2-bridge knots \cite{huh2021tight}.

Proven ropelength bounds typically apply in the limit of very large crossing number and are not competitive with measured upper bounds of simple knots from optimization algorithms. Proven lower bounds are typically several factors lower than the measured upper bounds. We recently showed that certain torus ropelengths do not follow asymptotic scaling relationships even at over 1000 crossings \cite{klotz2021ropelength}. Nevertheless, measured ropelength can help constrain conjectured bounds, and provide constraints on the coefficients of proven scaling relations.

Another feature of ideal knots is their space writhe, the writhe of the configuration averaged over every possible projection \cite{klenin2000computation}. The space writhe of ideal alternating knots was shown to be ``quasi-quantized'' near integer multiples of 4/7 \cite{pieranski2001quasi}, which follows from a ``predicted writhe'' invariant computed by a crossing nullification procedure on their diagrams \cite{cerf2000topological}. This suggests deeper connections between topological invariants and the geometric properties of ideal knots. By definition, topological invariants do not depend on the specific conformation of a knot, but the ideal configuration is unique (or part of a unique subset of configurations), and may be indicative of topological features of the knot. It was recently proven for example that the braid index of a knot may be used to establish a lower bound on the ropelength \cite{diao2020braid}, and a neural net evaluating the writhe representation of polymer knots was trained to predict their topology \cite{sleiman2023learning}. Correlations between knot invariants provide insight on underlying connections between them, and with a large enough sample of ideal knots we can investigate correlations between their geometric and topological features.

In this work, we use constrained gradient optimization to compute ropelength upper bounds for the 2176 knots with 12 minimal crossings, of which 1288 are alternating and 888 are non-alternating. Renderings of the ideal form of the biggest and smallest alternating and non-alternating knots are shown in Figure 1. Besides moving humanity's collective knowledge of ropelength from 11 crossings up to 12, this provides a large set of data with which to examine correlation between geometric and topological properties of knots. While most knots are non-alternating, most knots under 12 crossings are alternating. The investigation of writhe quantization \cite{pieranski2001quasi}, for example, used the 250 knots with under 10 crossings but only examined the alternating knots and did not have sufficient data to make inferences about non-alternating knots. With 888 12-crossing knots, this is no longer an issue.

Throughout the rest of this manuscript, we use the term ropelength to refer to the upper bounds on ropelength that we measured computationally. The ``true'' ropelength of all knots is unproven except for the unknot and certain classes of chained Hopf links \cite{cantarella2002minimum}. We also use the term ``writhe'' and ``space writhe'' to refer to the space writhe of the ideal configuration of the knot.

\begin{figure}
    \centering
    \includegraphics[width=0.5\textwidth]{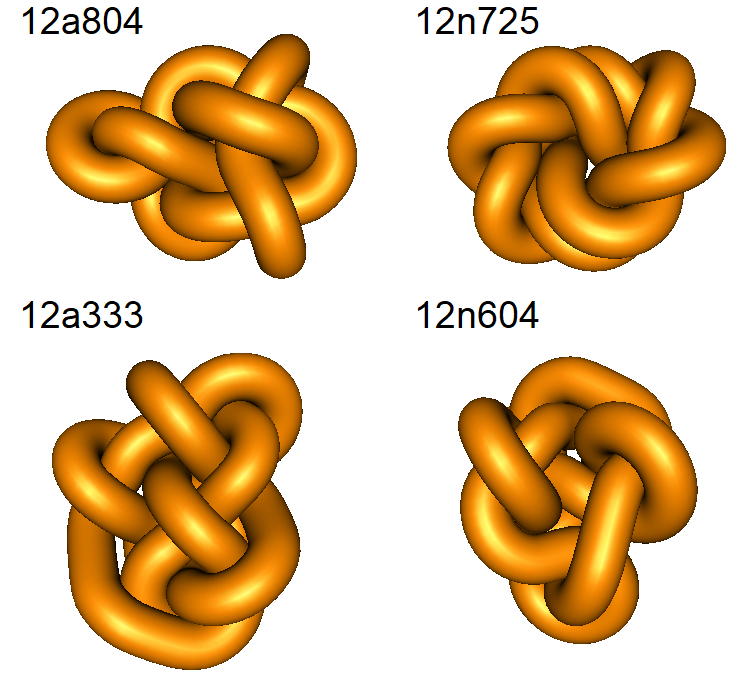}
    \caption{Configurations of the alternating (left) and non-alternating (right) knots with the smallest (top) and largest (bottom) ropelength, rendered using KnotPlot.}
    \label{fig:fourknots}
\end{figure}

\section{Methods}

We tighten knots using constrained gradient optimization implemented in the Ridgerunner software developed by Jason Cantarella \cite{ashton2011knot}. In short, it perturbs the coordinates of a discrete knot and follows the perturbations along a gradient towards minimal ropelength while respecting a no-overlap constraint. For initial configurations, we use Cartesian coordinates for 12-crossing knots that were generated by Se-Goo Kim and can be found on the KnotInfo database maintained by Charles Livingston and Allison Moore \cite{knotinfo}. The mean number of vertices of the configurations is 125 with a range between 88 and 160. These were downloaded and converted to Ridgerunner input files. Because it is not feasible to tune the run parameters to each individual knot, we used a common set of options in which knots were annealed with an equilateralization force for 25,000 steps or until the residual gradient reached 0.1, then the final configuration was run without the equilateralization force for 12,000 steps or until the residual reached 0.001. For knots with unfavorable initial conditions leading to deep local minima, we imported the coordinates into KnotPlot \cite{scharein2002interactive}, developed and maintained by Rob Scharein, and used a combination of Coulombic repulsion and a tangential tightening force to symmetrize and pre-tighten the configuration. We discuss the possibility of incomplete tightening and ropelength overestimation in the next section.



To identify anomalies, we compared the Alexander polynomials, evaluated at -1, of the initial and final configurations to each other and to their tabulated value. In six instances, the initial configurations were incorrect (and have since been updated on KnotInfo). In two instances, the knots underwent a topological change during the tightening process. These eight knots were re-sketched from diagrams and annealed using KnotPlot, before being re-entered into Ridgerunner. It remains a rare possibility that a knot is transformed into one with the same Alexander polynomial.




In addition to the ropelength, we calculated the space writhe and average crossing number from the ideal configurations. The space writhe is defined as the average writhe of the knot projected over every direction in space:

\begin{equation}
    Wr=\frac{1}{4\pi}\int_{K}\int_{K}\vec{dr_{1}}\times\vec{dr_{2}}\cdot\frac{\vec{r_{1}}-\vec{r_{2}}}{|\vec{r_{1}}-\vec{r_{2}} |^3},
\end{equation}
where $\vec{r_1}$ and $\vec{r_2}$ are the position vectors around the knot K. For discrete knots, it is computed as a double summation over every pair of line segments in the knot, as defined by Klenin and Langowski \cite{klenin2000computation}. The average crossing number is also defined over all possible projections of the knot and can be computed from: 

\begin{equation}
    ACN=\frac{1}{4\pi}\int_{K}\int_{K} \frac{|\hat{r}\cdot\vec{dr_{1}}\times\vec{dr_{2}}|}{|\vec{r_{1}}-\vec{r_{2}}|^2},
\end{equation}
where $\hat{r}$ is the unit vector along $\vec{r_2}-\vec{r_1}$. This was also computed as a double summation over a discrete knot. For computing space writhes of knots with 8 through 11 crossings, we used Brian Gilbert's Fourier coefficients to generate tight configurations with 512 vertices.



\section{Results and Discussion}

We found a population average for 12-crossing ropelengths of 102.95 (standard deviation 5.51), with the alternating knots having a mean of 106.989 (s.d 1.49) and the non-alternating 97.082 (s.d. 3.61). The mean ropelengths of each group of knots can be found in Table I at the end of this manuscript. A table of ropelengths, space writhes, and average crossing numbers may be downloaded from the Harvard Dataverse, as can Cartesian coordinates for each tight configuration \cite{dataverse}.  Histograms of 12-crossing ropelengths are shown in Figure 2, with the alternating knots having a tighter distribution. The alternating and non-alternating knots with the largest and smallest ropelength are shown in Figure 1. The parameters of the statistical distribution of ropelengths within a given crossing number and class are presented in Table II. Generally speaking, non-alternating knots have about 90\% the mean ropelength as alternating knots with the same crossing number, but at least twice the standard deviation. The average crossing number is weakly correlated with ropelength for alternating knots, and much more strongly for non-alternating knots (Pearson coefficients 0.30 and 0.81).

Figure 3 shows the ropelength distribution of knots between 3 and 12 crossings, with only the final distributions produced in this work. While these data lie visually on a line, we previously reported a best-fit power law of 0.81 ± 0.03, which has been updated to 0.83 ± 0.02. These fits are unweighted to avoid bias towards crossing numbers with more knots, and the uncertainties are standard errors on the fit parameters. The parameters from various fits to this data are presented in Table III. These fits should be treated as a heuristic way to estimate ropelength given a crossing number, but caution should be taken when making inferences about the ropelength scaling of the ensemble of all knots. 

While 12-crossing knots are typically not complex enough to interface with the proven lower bound for ropelength (which is currently 34.7 \cite{diao2003lower}), the knots studied in this work can in principle constrain the coefficients of power-law bounds. However, torus knots typically have the smallest ropelength in their class, and there are no torus knots with 12 crossings. For alternating knots, Diao recently proved that the ropelength must be greater than a linear lower bound with a slope of at least 1/59.5. The $11a367$ torus knot has a ropelength of at most 89.61, meaning that slope cannot exceed 8.15. However, the supercoiled configurations of Huh et al. \cite{huh2018ropelength} set a stronger bound of 7.63. It is known that non-alternating knots must have a ropelength greater than $1.1C^{3/4}$, with the true coefficient likely higher than 1.1. The $10_{124}$ torus knot provides an upper bound of 12.64 for the coefficient, slightly smaller than the bound provided by $12n725$ of 12.74 (although this may be reduced by annealing with more vertices). However, our previous work \cite{klotz2021ropelength} showed that the 1023 crossing T(32,33) knot constrains this to 10.76.

\begin{figure}
    \centering
    \includegraphics[width=0.7\textwidth]{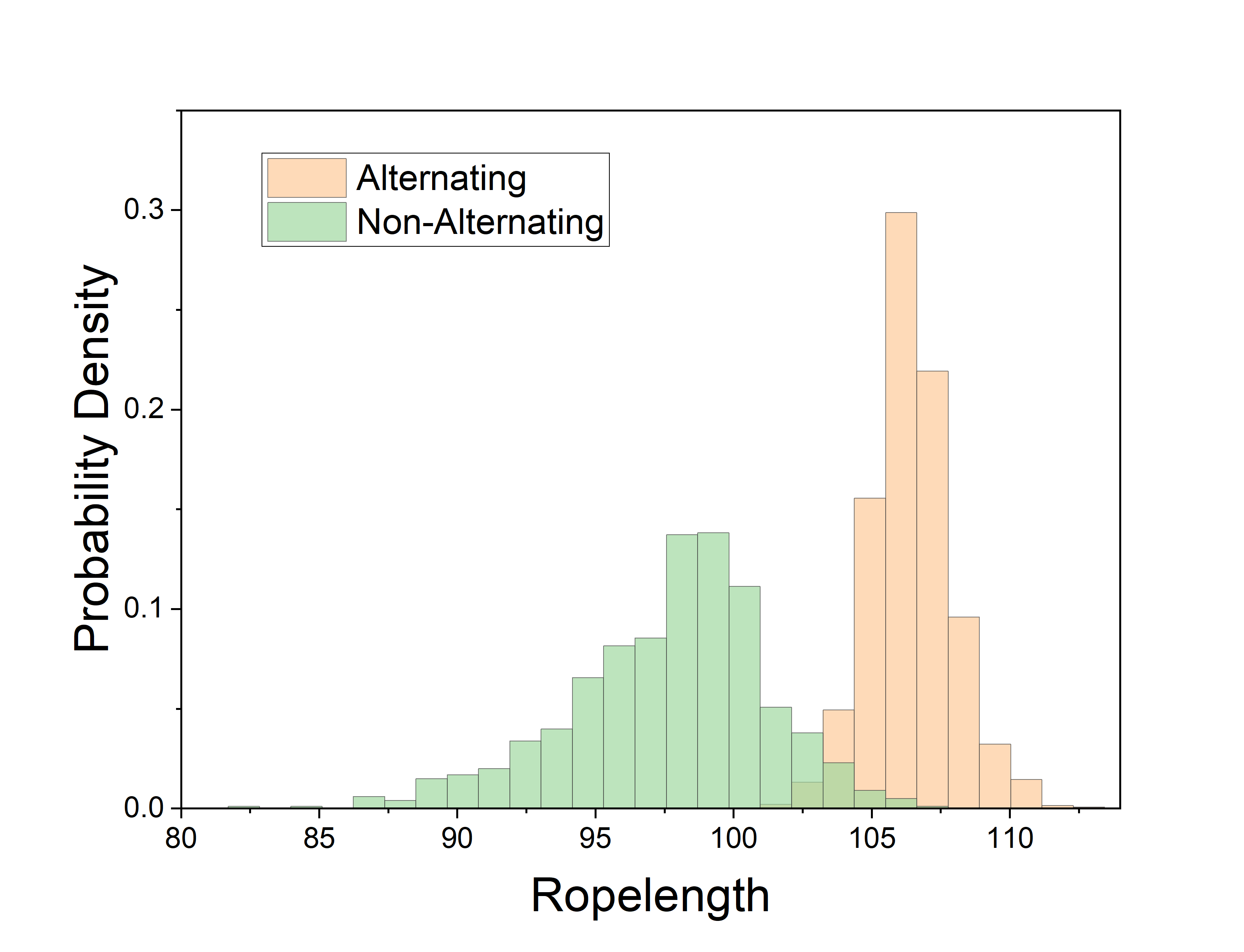}
    \caption{Normalized histograms of ropelengths for 12-crossing alternating and non-alternating knots. The non-alternating knots have a lower ropelength and a broader distribution, a feature seen for lower crossing numbers.}
    \label{fig:rlhist}
\end{figure}

\begin{figure}
    \centering
    \includegraphics[width=0.75\textwidth]{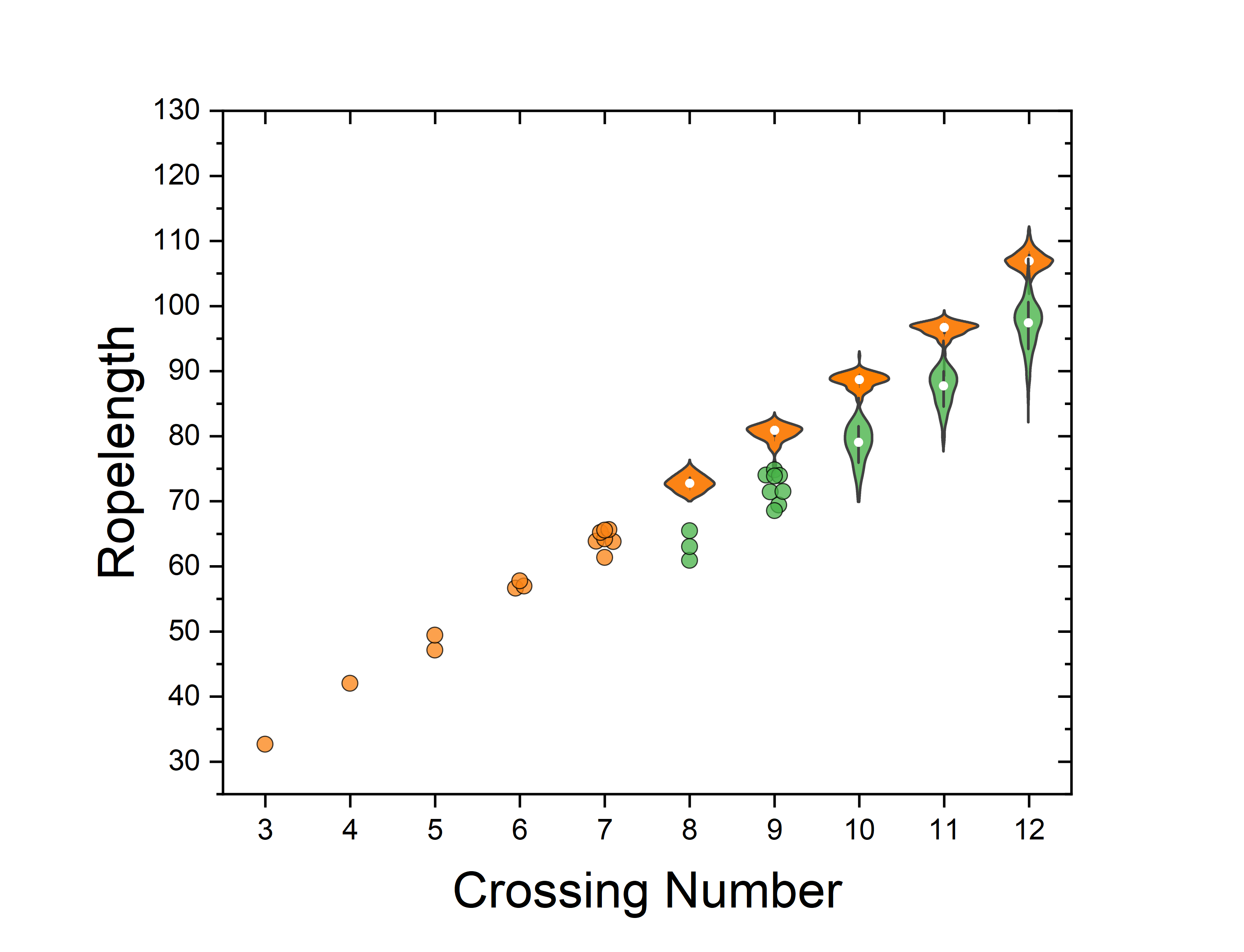}
    \caption{Plot of ropelengths as a function of crossing number, showing populations of alternating and non-alternating knots when appropriate. Crossing numbers with large populations are presented as violin distributions, while smaller populations are presented as scatter points, separated for visual clarity. Only the 12-crossing distributions are novel.}
    \label{fig:violin}
\end{figure}

\begin{figure}
    \centering
    \includegraphics[width=\textwidth]{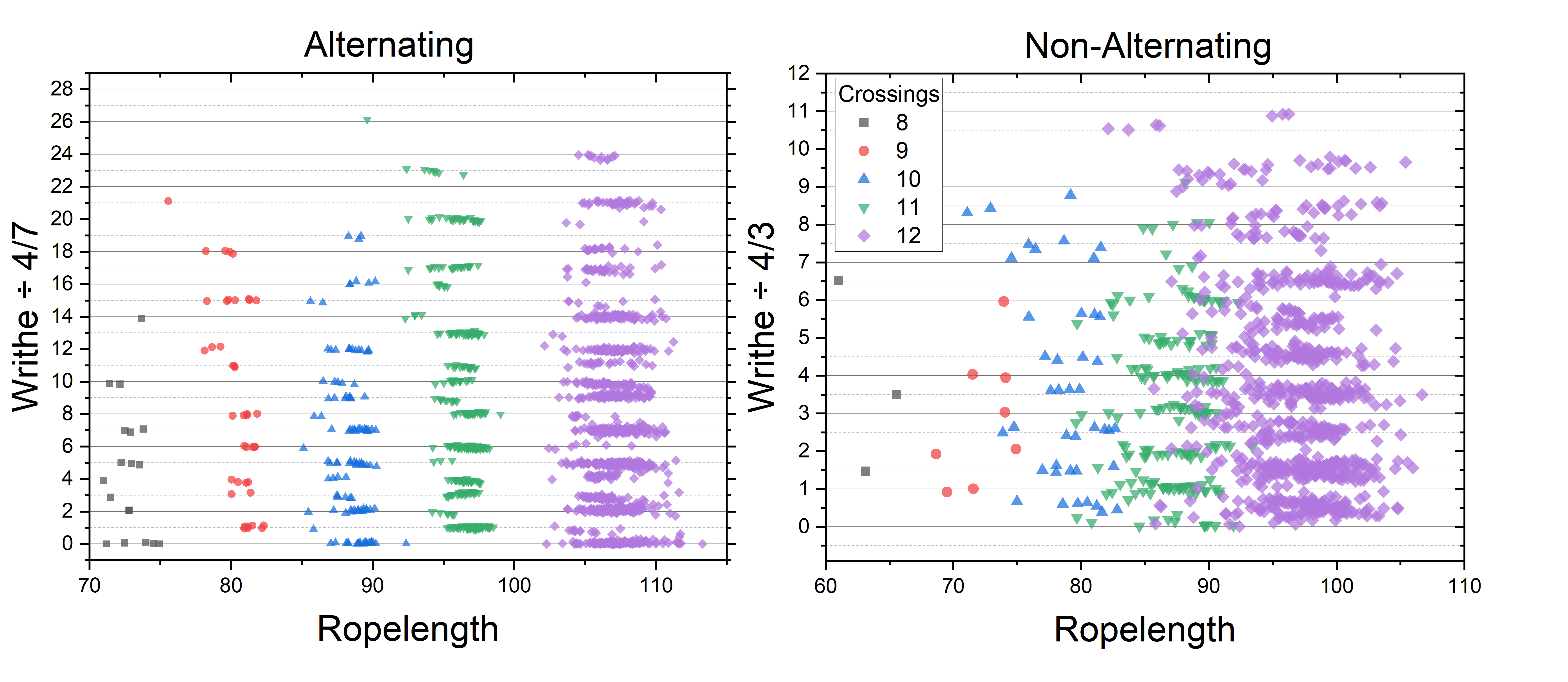}
    \caption{Scatter plots of space writhe against ropelength for alternating (left) and non-alternating (right) knots. The writhe has been scaled by its ``quantized'' value, 4/7 for alternating and 4/3 for for non-alternating, and is clustered at integer or half-integer values.}
    \label{fig:writhe}
\end{figure}

\begin{figure}
    \centering
    \includegraphics[width=0.7\textwidth]{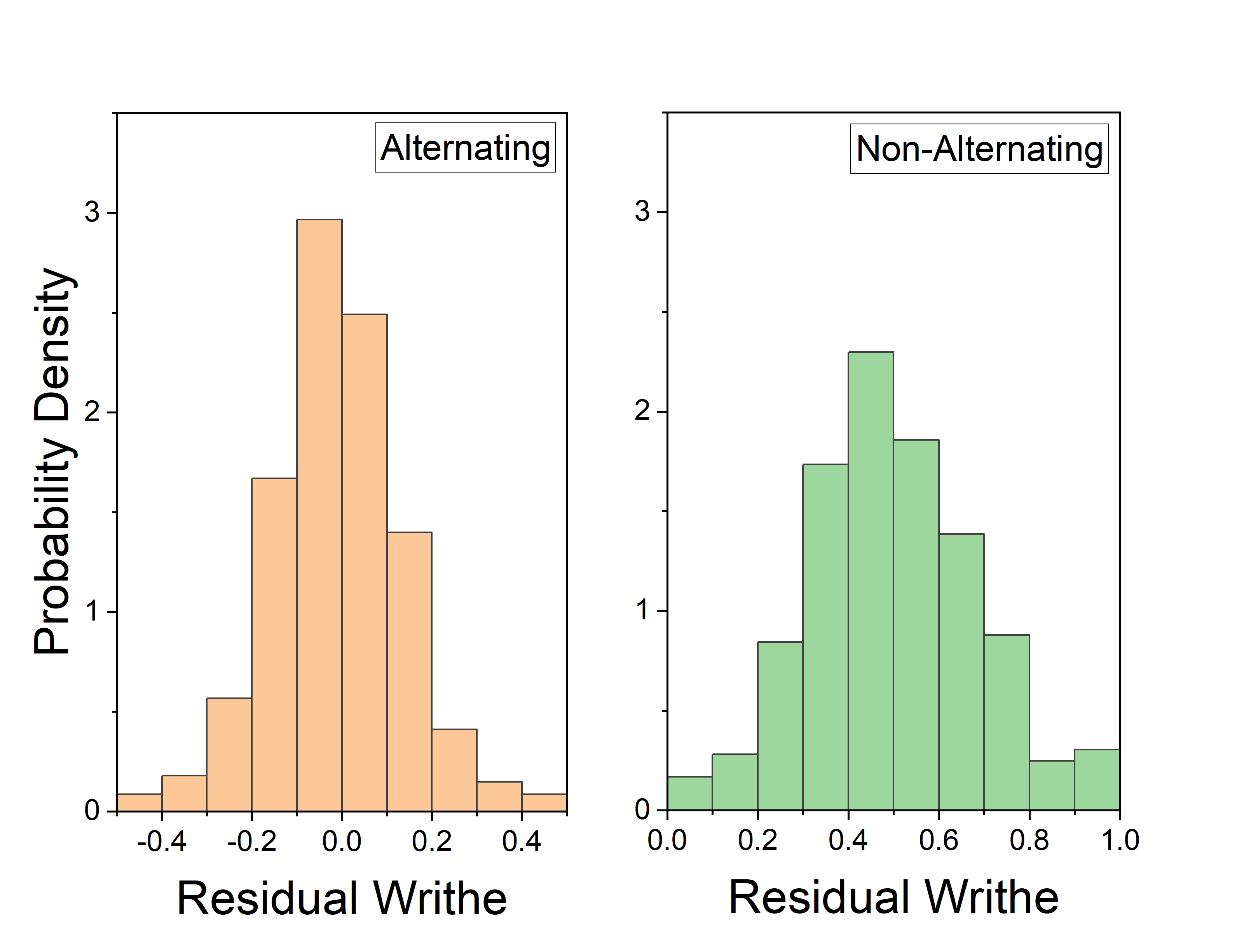}
    \caption{Normalized histograms of the deviation of the space writhe from its ``quantized'' value for 12-crossing alternating knots centered around 4/7 (left) and non-alternating knots centered around half-integer multiples of 4/3 (right).}
    \label{fig:writhehist}
\end{figure}

We computed the space writhe of our ideal 12-crossing knots. For ideal alternating knots, it was shown that reducing the diagram of a knot in a certain way produces an invariant, the predicted writhe, that predicts its ideal space writhe, which is an integer multiple of 4/7 \cite{cerf2000topological}. In a computational study, Pieranski and Pryzbwl verified that ideal alternating knots with up to 10 crossings do indeed have a writhe that is clustered near a multiple of 4/7. The majority of their tightened knots have a writhe that is within 0.1 of a multiple of 4/7, and the greatest deviation being 0.25. Only 20\% of knots with 10 or fewer crossings are non-alternating, and the writhe quantization of those knots was not investigated. The quantization of writhe is likely not absolute: the trefoil knot is predicted to have a writhe of 24/7, whereas the most precise measurement of its ideal form \cite{przybyl2014high} finds approximately 23.92/7. One should not necessarily assume that a knot with a larger residual writhe is not fully tight.


Our ideal alternating 12-crossing knots show evidence of clustering at multiples of 4/7, although there is more variance around the integer multiples than seen at lower crossings. Notably, we also observe clustering in the population of non-alternating knots when plot against their arbitrary index. A histogram of the space writhe shows peaks separated by a value greater than 1, without a peak at zero. To determine the periodicity of this clustering, we ran a two parameter sweep to minimize the variance of the fractional part of (Wr/A+B), finding that it is minimized when A=1.3434 and B=0.5253. We interpret this as having a quasi-quantization at half-integer multiples of 4/3, or odd multiples of 2/3. Applying this quantization scheme to non-alternating knots with 8 through 11 crossings, using ideal configurations generated by Gilbert for the Knot Atlas wiki, we found that knots with an even crossing number have a writhe near an odd multiple of 2/3, whereas knots with an odd crossing number have a writhe near an even multiple of 2/3. Within each crossing number, clusters of writhes are separated by 4/3. This is similar to the finding of Cerf and Stasiak \cite{cerf2000topological}, who found that even-component links have a predicted writhe that is a half-integer multiple of 4/7, with an integer multiple for an even component number.

The space writhe of the knots, scaled by their appropriate quantization, can be seen plot against their measured ropelength in Figure 4.  We do not assert a dependence between the two (they are uncorrelated), but use ropelength to space out the data for visualization. We have taken the absolute value of the space writhes for these plots, which is equivalent to choosing one stereoisomer for chiral knots. The alternating knots tend to be tightly clustered around their quantized values, whereas the non-alternating knots are more diffuse, but still clustered. The non-alternating knots with 8 and 9 crossings lie close to their quantized values, while deviations grow for 10 through 12 crossings. A crossing nullification procedure was used by Stasiak to show that the $8_{20}$ knot has a predicted writhe of 2 \cite{stasiak2000knots}, which is a half-integer multiple of 4/3. Three independent measurements \cite{katlas,ashton2011knot,pieranski1998search} find an ideal writhe of 1.95 or 1.96.


The distributions of the ideal space writhes of 12-crossing knots about their quantized values are shown in Figure 5. Alternating knots about 4/7 have a mean of -.01 and a standard deviation of 0.14, non-alternating knots about 4/3 have a mean of .50 and a standard deviation of 0.19. For comparison, the Pieranski and Pryzbwl \cite{pieranski2001quasi} measured a standard deviation of approximately 0.09 for alternating knots with 10 or fewer crossings. The null hypothesis of a uniform distribution between 0 and 1 has a standard deviation of 0.29. It is possible that there may be exceptions to the non-alternating half-integer rule, as indicated by the top three points of Figure 4 and the rightmost bin of Figure 5. Knots with large residual writhe may be candidates for enhanced shrinking to eliminate outliers, but this did not prove a fruitful strategy for reducing the mean measured ropelength.

It is thought that correlations between knot invariants are indicative of deeper connections between properties of knots, particular between algebraic and geometric invariants \cite{davies2021advancing}. For example, the Volume Conjecture posits a direct relationship between the hyperbolic volume of a knot and its Jones polynomial \cite{murakami2018volume}. An ideal knot has two uncorrelated properties, ropelength and writhe, that may be correlated with other topological invariants of the knot. Correlations between a specific conformation of a knot and its invariants may be counterintuitive, as a knot may be swollen to any length and twisted to have an arbitrary writhe, and one may not expect, for example, that a polynomial derived from a diagram of a knot ``knows'' how tight it can become. With a sample of over 2000 knots within a single crossing number, we can examine correlations between ropelength and space writhe and other knot invariants. We find that the ropelength is correlated with hyperbolic volume (Figure 6), with Pearson coefficients of 0.35 and 0.74 for alternating and non-alternating knots. The ropelength is also correlated with the -1 evaluation of the Alexander polynomial (Pearson coefficients 0.39 and 0.55), although that is already known to be correlated with hyperbolic volume \cite{friedl2011approximations}. We examine the same relations in lower-crossing knots, we find very strong correlations that had not been examined previously, for example 0.95 for 10-crossing alternating knots. In a sense, a knot's hyperbolic volume is a measure of its complexity determined by the number of tetrahedra that must be glued together to construct its complement, and the ropelength is a measure of complexity determined by the number of curved strings that must be glued together, so a correlation is not outside the realm of possibility. The reduced correlation with increased crossing number may be indicative that this trend is a property of simple knots and vanishes as knots become asymptotically complex, or just an indication that it is more difficult to converge on the ropelength of more complex knots.

While all 12-crossing knots are hyperbolic, it is conceivable that ropelength and the Alexander polynomial are correlated for torus knots. For example, alternating torus knots have both a ropelength and an Alexander(-1) that is linear in crossing number. The lowest crossing number that admits three torus knots is 63, and we verified that the ropelength and absolute value of Alexander(-1) are ordered \footnote{Specifically, T(9,8) has an Alexander value of 9 and a ropelength of 260.0, T(21,4) has an Alexander of 21 and a ropelength of of 275.5, and T(63,2) has an Alexander of -63 and a ropelength likely between 486 and 532 \cite{huh2018ropelength}}. To examine correlations within a single crossing number for torus knots, one would likely have to examine knots with 4320 crossings, of which 10 are torus knots \cite{oeis}. We make the additional observation that the torus knots with 5 through 11 crossings each have the smallest ropelength in their class, and zero hyperbolic volume.

\begin{figure}
    \centering
    \includegraphics[width=\textwidth]{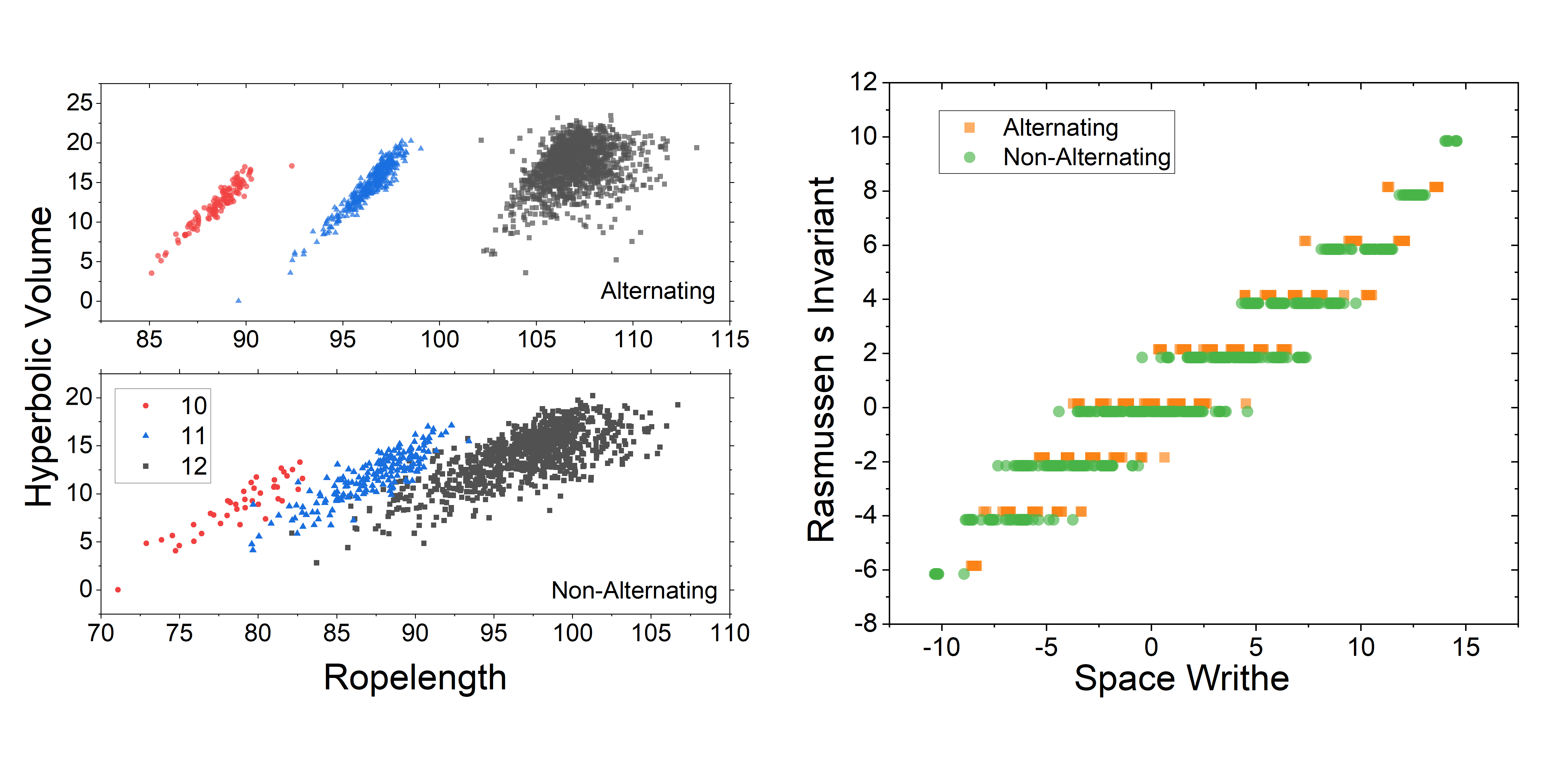}
    \caption{Correlation between geometric and topological properties of tight knots. Left: Hyperbolic volume plot against ropelength for 10, 11, and 12-crossing knots showing positive correlation. Right: Rasmussen s invariant plot against space writhe for 12-crossing knots, showing positive correlation. Points are shifted vertically for visual clarity.}
    \label{fig:invariants}
\end{figure}

Since ideal writhe is uncorrelated with ropelength, we may independently examine its correlation with other invariants. We find that the ideal writhe shows a strong correlation (0.96) with the Rasmussen s and Ozsvath-Szabo tau invariants, the former shown in Figure 6. The correlation is present for 11-crossing knots, but does not appear significant for 10. These parameters arise from Khovanov homology and bound the 4-genus of the knots. We have no formal or heuristic explanation for this.

\subsection{Looseness Considerations}

There are several heuristic reasons to suspect we may be overestimating ropelengths. The 12-crossing point lies above the trendline in all the fits in Table III. The 12-crossing alternating knots are the only large sample to have a positive skew. There is a weak positive correlation between the alternating ropelength and its residual writhe, and the non-alternating writhe is considerably scattered about its quanta (however, even extremely tight knots are known not to lie exactly on a writhe quantum). There are three main reasons that our Ridgerunner tightenings would not reach a value close to the minimum possible: insufficient vertices, insufficient runtime, and unfavorable initial conditions.

The mean number of vertices used in our measurements is 125 but they range between 88 and 160. To ascertain the effect of vertex number, we can consider a 12-crossing Hopf chain-link necklace consisting of six closed loops, which has a known ropelength of $6\cdot(4\pi+4)\approx99.4$ \cite{cantarella2002minimum}. With 20 vertices per link (120 total), the same Ridgerunner parameters shrink it to 100.63, 1.2\% above the true value. With 12 vertices per link, Ridgerunner reaches 102.24, 2.9\% above the minimum. Based on vertex count alone, it would be safe to estimate at least 1\% ropelength excess. A negative correlation between vertex number and ropelength may indicate a bias, but there is no correlation for alternating knots and a weak positive correlation for non-alternating knots.

We do not believe insufficient runtime has a large effect on our population ropelength. Tightenings that terminate at the pre-set residual did not have a significantly smaller ropelength than those that terminated after 12,000 steps. We ran several of the knots with the largest ropelength for an extra 2000 steps, and found that only the second or third decimal digit was affected.

Unfavorable initial conditions can cause Ridgerunner to get stuck in deep local minima, and this is likely the biggest contributor to the excess ropelength in our measurements. After our initial run, several knots had an outlying ropelength, far beyond the expected tail of the histograms. This is likely due to a wide range of distances between adjacent strands in the knots, which causes Ridgerunner to tighten regions that are already tight while leaving others excessively loose. The initial coordinates of these outliers were imported into KnotPlot, symmetrized using Coulombic repulsion, pre-tightened, and then re-entered into Ridgerunner. Doing this for obvious outliers reduced the population average by about 0.5. Readers who wish to continue this process may find our ropelengths and coordinates on the Harvard Dataverse \cite{dataverse}.

\section{Conclusion}

We have measured ropelength upper bounds and ideal writhes for the 2176 12-crossing knots. Although these ropelength measurements are likely at least 1\% overestimated, they add to the data available for understanding the relationship between ideal knot configurations and topological invariants. Notably, we find evidence that non-alternating knots are clustered around integer and half-integer values of 4/3, that ropelength is correlated with hyperbolic volume, and that ideal writhe is correlated with the Rasmussen s invariant from Khovanov homology.

There are 9,988 knots with 13 crossings, which is the first crossing number at which most knots are non-alternating. Unless the algorithms or available hardware are significantly improved, we do not recommend an exhaustive ropelength for higher crossing numbers. However, it would be worthwhile to investigate a subset of a few hundred 13 or higher-crossing knots to determine whether the trends observed in this work continue as knot complexity increases. Readers who wish to undertake such a task should be advised that an approach combining Coulomb repulsive molecular dynamics and a traditional minimization algorithm may be more efficient than simply shrinking knots from their initial coordinates.

\section{Acknowledgements}
This work is supported by the National Science Foundation, grant number 2105113. We thank Andrzej Stasiak for helpful discussions, and Allison Moore, Charles Livingston, and Se-Goo Kim for providing the knot coordinates. We are generally grateful to everyone who develops, maintains, and answers emails about freely available tools for knot theory investigations, including Rob Scharein, Luca Tubiana, Dror Bar-Natan, Jason Cantarella, Brian Gilbert and others.

\section{Tables}

\begin{table}[h]
\begin{tabular}{lllllll}
\hline
\multicolumn{1}{|l|}{Crossings} & \multicolumn{1}{l|}{Population} & \multicolumn{1}{l|}{s.d} & \multicolumn{1}{l|}{Alt.} & \multicolumn{1}{l|}{s.d} & \multicolumn{1}{l|}{Non-Alt.} & \multicolumn{1}{l|}{s.d} \\ \hline
\multicolumn{1}{|l|}{3} & \multicolumn{1}{l|}{32.74} & \multicolumn{1}{l|}{} & \multicolumn{1}{l|}{32.74} & \multicolumn{1}{l|}{} & \multicolumn{1}{l|}{} & \multicolumn{1}{l|}{} \\ \hline
\multicolumn{1}{|l|}{4} & \multicolumn{1}{l|}{42.09} & \multicolumn{1}{l|}{} & \multicolumn{1}{l|}{42.09} & \multicolumn{1}{l|}{} & \multicolumn{1}{l|}{} & \multicolumn{1}{l|}{} \\ \hline
\multicolumn{1}{|l|}{5} & \multicolumn{1}{l|}{48.32} & \multicolumn{1}{l|}{1.6} & \multicolumn{1}{l|}{48.32} & \multicolumn{1}{l|}{1.6} & \multicolumn{1}{l|}{} & \multicolumn{1}{l|}{} \\ \hline
\multicolumn{1}{|l|}{6} & \multicolumn{1}{l|}{57.16} & \multicolumn{1}{l|}{0.58} & \multicolumn{1}{l|}{57.16} & \multicolumn{1}{l|}{0.58} & \multicolumn{1}{l|}{} & \multicolumn{1}{l|}{} \\ \hline
\multicolumn{1}{|l|}{7} & \multicolumn{1}{l|}{63.37} & \multicolumn{1}{l|}{2.85} & \multicolumn{1}{l|}{63.37} & \multicolumn{1}{l|}{2.85} & \multicolumn{1}{l|}{} & \multicolumn{1}{l|}{} \\ \hline
\multicolumn{1}{|l|}{8} & \multicolumn{1}{l|}{71.64} & \multicolumn{1}{l|}{3.85} & \multicolumn{1}{l|}{72.78} & \multicolumn{1}{l|}{1.13} & \multicolumn{1}{l|}{63.21} & \multicolumn{1}{l|}{2.27} \\ \hline
\multicolumn{1}{|l|}{9} & \multicolumn{1}{l|}{79.30} & \multicolumn{1}{l|}{4.06} & \multicolumn{1}{l|}{80.45} & \multicolumn{1}{l|}{1.31} & \multicolumn{1}{l|}{72.26} & \multicolumn{1}{l|}{2.33} \\ \hline
\multicolumn{1}{|l|}{10} & \multicolumn{1}{l|}{85.98} & \multicolumn{1}{l|}{4.64} & \multicolumn{1}{l|}{88.54} & \multicolumn{1}{l|}{1.16} & \multicolumn{1}{l|}{78.81} & \multicolumn{1}{l|}{2.81} \\ \hline
\multicolumn{1}{|l|}{11} & \multicolumn{1}{l|}{93.38} & \multicolumn{1}{l|}{4.66} & \multicolumn{1}{l|}{96.43} & \multicolumn{1}{l|}{1.15} & \multicolumn{1}{l|}{87.35} & \multicolumn{1}{l|}{2.71} \\ \hline
\multicolumn{1}{|l|}{12} & \multicolumn{1}{l|}{102.95} & \multicolumn{1}{l|}{5.51} & \multicolumn{1}{l|}{106.99} & \multicolumn{1}{l|}{1.49} & \multicolumn{1}{l|}{97.08} & \multicolumn{1}{l|}{3.61} \\ \hline
 &  &  &  &  &  & 
\end{tabular}
\caption{Mean ropelengths of knots between 3 and 12 crossings, with standard deviations.}
\end{table}

\begin{table}[h]
\begin{tabular}{|l|l|l|l|}
\hline
Crossings & s.d/mean & Skewness & Kurtosis \\ \hline
8A & 0.0155 & 0.1273 & 2.2034 \\ \hline
9A & 0.0162 & -1.5173 & 6.1508 \\ \hline
10A & 0.0131 & -0.4959 & 3.7989 \\ \hline
10N & 0.0357 & -0.7594 & 3.0249 \\ \hline
11A & 0.0119 & -1.3971 & 6.9911 \\ \hline
11N & 0.0311 & -0.6620 & 3.0830 \\ \hline
12A & 0.0139 & 0.2559 & 3.7451 \\ \hline
12N & 0.0372 & -0.5964 & 3.6611 \\ \hline
\end{tabular}
\caption{Properties of ropelength distributions within each crossing number, for samples with 18 or more knots.}
\end{table}

\begin{table}[h]
\begin{tabular}{|l|l|l|}
\hline
Sample & Linear & Power \\ \hline
All Knots &  $(7.6\pm0.1)C+(10.7\pm0.7)$ & $(12.9\pm0.4)C^{0.83\pm0.02}$ \\ \hline
Alternating & $(8.1\pm0.1)C+(8.5\pm0.9)$ & $(12.0\pm0.5)C^{0.87\pm0.02}$ \\ \hline
Non-Alternating & $(8.3\pm0.3)C-(3.1\pm3.0)$ & $(7.2\pm0.6)C^{1.04\pm0.04}$ \\ \hline
\end{tabular}
\caption{Linear and power fits to the ropelength data.}
\end{table}
\color{white}
I\\

\color{black}
\bibliographystyle{unsrt}
\bibliography{knotrefs2}
\end{document}